\documentclass[runningheads]{llncs}
\pdfoutput=1

\usepackage[T1]{fontenc}
\usepackage{graphicx}

\usepackage{color}
\usepackage[colorlinks]{hyperref}

\usepackage{cite}

\usepackage{amsfonts}
\usepackage{mathtools}
\usepackage{booktabs}
\usepackage[textsize=footnotesize]{todonotes}
\usepackage{siunitx}

\DeclarePairedDelimiterX{\set}[1]\lbrace\rbrace{\setaux #1||\endsetaux}
\def\setaux#1|#2|#3\endsetaux{\if\relax\detokenize{#2}\relax #1 \else #1 \;\delimsize\vert\; #2 \fi}

\newcommand{\integers}{\ensuremath{\mathbb Z}}

\newcommand{\B}[1]{\{0,1\}^{#1}}

\newcommand{\X}{\mathcal{X}}
\newcommand{\bigO}{\ensuremath{\mathcal O}}

\newcommand{\Sym}[1]{\ensuremath{\mathrm{Sym}_{#1}}}

\newcommand{\T}{^\top}
\newcommand{\sprod}[2]{{#1}\T{#2}}
\newcommand{\define}{\coloneqq}

\begin{document}
    \title{Handling Sub-symmetry in Integer Programming using Activation Handlers}
    %
    %
    \author{Christopher Hojny
      \and
      Tom Verhoeff
      \and
      Sten Wessel
    }
    \authorrunning{C.\ Hojny et al.}
    %
    \institute{Eindhoven University of Technology, Eindhoven, The Netherlands\\\email{\{c.hojny,t.verhoeff,s.wessel\}@tue.nl}}

    \maketitle

    \begin{abstract}
        Symmetry in integer programs (IPs) can be exploited in order to reduce solving times.
        Usually only symmetries of the original IP are handled, but new symmetries may arise at some nodes of the branch-and-bound tree.
        While symmetry-handling inequalities (SHIs) can easily be used to handle original symmetries, handling sub-symmetries arising later on is more intricate.
        To handle sub-symmetries, it has recently been proposed to add SHIs that are activated by auxiliary variables.
        This, however, may increase the size of the IP substantially as all sub-symmetries need to be modeled explicitly.
        As an alternative, we propose a new framework for generically activating SHIs, so-called \emph{activation handlers}.
        This framework allows for a direct implementation of routines that check for active sub-symmetries, eliminating the need for auxiliary variables.
        In particular, activation handlers can activate symmetry-handling techniques that are more powerful than SHIs.
        We show that our approach is flexible, with applications in the multiple-knapsack, unit commitment, and graph coloring problems.
        Numerical results show a substantial performance improvement on the existing sub-symmetry-handling methods.

    \end{abstract}

    \section{Introduction} \label{sec:introduction}

    One of the most popular methods to solve integer programs is
    branch-and-bound (B\&B), which iteratively splits the integer program into smaller
    subproblems that are solved in turn~\cite{LandDoig1960}.
    While B\&B can solve problems with thousands of variables and
    constraints in adequate time, it is well-known that the presence of
    symmetries leads to unnecessarily large B\&B trees.
    The main reason is that B\&B explores symmetric
    subproblems, which all provide essentially the same information.
    Therefore, symmetry-handling is an important ingredient of modern B\&B
    implementations that substantially improves the running
    time~\cite{pfetsch2019computational}.
    
    We consider permutation symmetries of binary
    programs~$\max\set{ \sprod{c}{x} | x \in \X }$, where~$c \in
    \integers^d$ and~$\X \subseteq \B{d}$.
    A \emph{permutation} is a bijection of~$[d] \define \set{1,\dots,d}$;
    the set of all permutations is~$\Sym{d}$.
    We assume that a permutation~$\pi \in \Sym{d}$ acts on~$x \in
    \B{d}$ by permuting its coordinates, i.e., $\pi(x) \define
    (x_{\pi^{-1}(1)},\dots,x_{\pi^{-1}(d)})$.
    A \emph{symmetry} of the binary program is a permutation~$\pi \in
    \Sym{d}$ that preserves
    the objective, i.e., $\sprod{c}{\pi(x)} = \sprod{c}{x}$, and
    feasibility, i.e., $x \in \X$ if and only if~$\pi(x) \in \X$.
    Note that the set of all symmetries forms a group under composition.
    We refer to this group as the \emph{symmetry group} of the binary
    program, denoted by~$G$.
    Since computing~$G$ is NP-hard~\cite{margot2010symmetry}, one usually only handles a
    subgroup of~$G$, which can either be detected
    automatically~\cite{pfetsch2019computational,salvagnin2005dominance} or
    is provided by an expert.
    
    To solve binary programs, B\&B generates
    subproblems~$\max\set{\sprod{c}{x} | x \in Q}$, where~$Q \subseteq
    \X$.
    Each subproblem is a binary program with symmetry
    group~$G_Q$.
    In general, $G_Q$ is different from~$G$ and neither is a
    subgroup of the other.
    Following~\cite{bendotti2020symmetry}, we call symmetries in~$G_Q$
    \emph{sub-symmetries} of the initial binary program.
    If sub-symmetries appear frequently during the solving process, it can
    be beneficial to handle them.
    But since computing (subgroups of)~$G_Q$ might be costly,
    providing knowledge about sub-symmetries via experts can substantially
    reduce the complexity of handling sub-symmetries.
    This, however, leads to a new challenge: how to efficiently provide
    expert knowledge to a binary programming solver.

    Recently, \cite{bendotti2020symmetry} suggested to introduce, for each possible
    sub-symmetry, a simple symmetry-handling inequality  (SHI) that is
    coupled with auxiliary variables that enable (resp.\
    disable) the SHI if the sub-symmetry is active (resp.\ inactive) at a
    subproblem.
    This, however, might lead to a significant increase in the size of the problem
    formulation as all sub-symmetries need to be explicitly modeled.
    Moreover, simple SHIs might not lead to the strongest symmetry
    reductions.

    In this paper, we propose an alternative approach.
    Instead of using auxiliary variables to activate SHIs, we introduce a
    simple framework that allows to activate arbitrary symmetry-handling
    techniques.
    Our so-called \emph{activation handlers} receive expert knowledge as
    input, i.e., rules to decide which sub-symmetries are present at a
    subproblem.
    These rules are then automatically evaluated and communicate with the
    solver how sub-symmetries are handled.
    We thus mimic the process of how an expert would handle sub-symmetries and
    avoid reformulating the original problem.
    In particular, our framework is more flexible than auxiliary variables
    as not all activation rules can be compactly expressed by variables. 

    In the remainder of this section, we provide basic notation as
    well as terminology and we provide a brief overview of symmetry
    handling methods.
    Sec.~\ref{sec:subsym} summarizes the state-of-the-art of
    handling sub-symmetries, which is complemented by a description of our
    activation handler framework in Sec.~\ref{sec:ah}.
    Then, we illustrate for three classes of problems how activation
    handlers can be used to handle sub-symmetries
    (Sec.~\ref{sec:applications}),
    and we evaluate our activation handler framework on a broad set of
    instances (Sec.~\ref{sec:experiments}).
    Numerical results show that our novel framework substantially improves
    upon the state-of-the-art in handling sub-symmetries.

    \paragraph{Notation and Terminology}
    Let~$G$ be the symmetry group of a binary program with feasible
    region~$\X$.
    The \emph{orbit} of~$x$ is~$\mathrm{orb}_G(x) = \set{\pi(x) | \pi \in G}$ and contains all solutions equivalent to~$x$ w.r.t.~$G$.
    Note that the orbits partition~$\mathcal{X}$.

    To handle symmetries, it suffices to restrict~$\X$ to
    (usally lexicographically maximal)
    representatives of orbits.
    A vector~$y \in \integers^d$ is \emph{lexicographically greater than} $z \in \integers^d$, denoted~$y \succ z$, if there is~$i \in [d]$ such that $y_i > z_i$ and $y_j = z_j$~for all~$j < i$.
    We write~$y \succeq z$ when~$y \succ z$ or~$y = z$ hold.
    Then, a solution~$x \in \mathcal X$ is \emph{lexicographically maximal} in its
    orbit under~$G$ when~$x \succeq \pi(x)$ for all~$\pi \in G$.

    \paragraph{Related Literature}
    A great variety of symmetry-handling methods exist in the literature,
    including variable branching and fixing
    rules~\cite{margot2003exploiting,OstrowskiEtAl2011}, pruning
    rules~\cite{Margot2002,margot2003exploiting,ostrowski2009symmetry},
    model reformulation techniques~\cite{FischettiLiberti2012}, and symmetry-handling
    constraints~\cite{Friedman2007,hojny2019polytopes,Hojny2020,kaibel2008packing,Liberti2008,Liberti2012,Liberti2012a,LibertiOstrowski2014,Salvagnin2018}.
    In the following, we provide details about some constraint-based
    techniques that we also use in our experiments.
    For details on other techniques, we refer to the
    overview~\cite{margot2010symmetry} and the computational
    survey~\cite{pfetsch2019computational}.

    The common ground of most symmetry-handling constraints is to enforce
    that only lexicographically maximal representatives of symmetric
    solutions are computed.
    For a single symmetry~$\pi \in G$, $x \succeq \pi(x)$ can
    be enforced in linear time by propagation and separation
    techniques~\cite{DoornmalenHojny2022,hojny2019polytopes}.
    For general groups, however, it is coNP-complete to decide whether a
    solution is lexicographically maximal in its orbit if~$G$ is given by a
    set of generators~\cite{babai1983canonical}.
    Attention has thus been spent on groups that arise frequently in
    practice.
    One such case assumes that the variables are organized in a matrix~$x =
    (x_{i,j})_{i \in [m],j \in [n]}$ and that the symmetries in~$G$ permute
    the columns of~$x$ arbitrarily.
    Such symmetries arise frequently in benchmark
    instances~\cite{pfetsch2019computational},
    and in Sec.~\ref{sec:applications}, we illustrate some applications.
    There, we also discuss how orbitopes can be used to handle
    these symmetries.

    The \emph{full orbitope} is the convex hull of all binary matrices
    with lexicographically non-increasingly sorted columns.
    If one restricts to matrices all of whose rows have at most
    (resp.\ exactly) one~1-entry, the corresponding convex hull is called
    \emph{packing orbitope} (resp.\ \emph{partitioning orbitope}).
    The above-mentioned matrix symmetries can be handled by separating
    valid inequalities for orbitopes.
    For packing/partitioning orbitopes, a facet description can be
    separated in linear time~\cite{kaibel2008packing}; for full orbitopes,
    efficiently separable IP formulations are
    known~\cite{hojny2019polytopes}.
    Moreover, efficient propagation algorithms for both full and
    packing/partitioning orbitopes are
    known~\cite{bendotti2021orbitopal,KaibelEtAl2011}.
    The so-called \emph{orbitopal fixing} algorithms receive local variable
    bounds at a subproblem and derive variables that need to be fixed at a
    certain value to guarantee that a solution is contained in the
    orbitope.

    \section{Sub-symmetry in Integer Programming} \label{sec:subsym}

    In this section, we provide a detailed explanation of sub-symmetries.
    We discuss how sub-symmetries arise at a
    subproblem~$\max\set{c\T x | x \in Q}$ (Sec.~\ref{sec:subsymbinary}) and how one can handle them
    by state-of-the-art techniques (Sec.~\ref{sec:handlesubsym}).
    We use the multiple knapsack problem as a running example for
    illustration purposes.

    \begin{example}\label{ex:mkp}
        The \emph{multiple knapsack problem}~(MKP) considers~$m$~items with associated profit~$p_i$ and weight~$w_i$, $i \in [m]$, as well as~$n$~knapsacks with capacity~$c_j$, $j \in [n]$.
        The objective is to assign each item to at most one knapsack such that the total weight of items assigned to knapsack~$j \in [n]$ does not exceed~$c_j$ and the total profit of assigned items is maximized.
        The MKP is NP-hard as it generalizes the NP-hard single-knapsack problem~\cite{GareyJohnson1979}.
        A standard IP formulation is given by~\ref{eq:ip:mks}, which is given below.
        There,variable~$y_{i,j}$ indicates whether item~$i \in [m]$ is packed in knapsack~$j \in [n]$, see~\cite{MARTELLO1987213}.
        Let~$\mathcal Y_\text{MKP}$ denote the set of all feasible solution matrices~$y$.

        This formulation exhibits multiple types of symmetries.
        From any feasible solution, permuting indices of knapsacks with equal capacity yields another feasible solution.
        For example, when all knapsacks have the same capacity~$c = c_j$ for all~$j \in [n]$, this symmetry corresponds to permutations the columns of the solution matrix~$y$.
        Symmetries also arise from items with identical properties, i.e., identical weight and profit.
        In any feasible solution these items can be permuted, corresponding to permutations of the respective rows in the solution matrix~$y$.

        Sub-symmetries also occur in this formulation, arising from the following observation on partially-filled knapsacks.
        Consider two knapsacks~$j$ and $j'$, and an item with index~$i$.
        Suppose that items~$\set{1,\dots,i-1}$ are placed such that the capacity that remains in knapsacks~$j$ and~$j'$ are equal.
        Then, the placement of the remaining items~$\set{i, \dots, m}$ can be permuted between the two knapsacks.
        We call this type of sub-symmetry \emph{capacity sub-symmetries}.
        Note that by this definition also knapsacks~$j$ and~$j'$ with~$c_j \neq c_{j'}$ can become sub-symmetric.
        \begin{subequations} \label{eq:ip:mks}
            \begin{alignat}{5}
                \text{maximize}\quad && \sum_{i=1}^m \sum_{j=1}^n p_i y_{i,j} \span \label{eq:ip:mks:start} \\
                \text{subject to}\quad && \sum_{i=1}^m w_i y_{i,j} &\le c_j &&\quad\forall_{j \in [n]}, \label{eq:ip:mks:capacity} \\
                && \sum_{j=1}^n y_{i,j} &\le 1 &&\quad\forall_{i \in \mathcal [m]}, \label{eq:ip:mks:packing} \\
                && y_{i,j} &\in \set{0, 1} &&\quad\forall_{i \in [m], j \in [n]}, \label{eq:ip:mks:end}
              \end{alignat}
        \end{subequations}
    \end{example}

    \subsection{Sub-symmetries in General Binary Programs}\label{sec:subsymbinary}

    In general, sub-symmetries can be defined for arbitrary collections of subproblems~$\mathbb S = \set{Q_s \subset \mathcal X | s \in [q]}$ for some~$q$.
    The sub-symmetries then correspond to permutations in~$G_{Q_s}$, $s \in [q]$, and are either automatically detected by an IP solver or are provided by a user.
    In the example of the MKP, the solution subsets for which the capacity sub-symmetries occur can be defined as
    \begin{equation}
        Q^i_{j,j'} = \left\{y \in \mathcal Y_\text{MKP} \;\middle|\; c_j - \sum_{k=1}^{i-1} w_k y_{k,j} = c_{j'} - \sum_{k=1}^{i-1} w_k y_{k,j'}\right\}
    \end{equation}
    for all pairs~$j,j' \in [n]$, $j < j'$, and all~$i \in [m]$.
    We denote by~$\mathbb S_\text{MKP}$ the collection of all these solution subsets.

    During the branch-and-bound process, a subproblem corresponding to a node in the B\&B~tree may belong to one or multiple solution subsets~$Q_s$.
    We then say that the sub-symmetries in~$G_{Q_s}$ become \emph{active}.
    To exploit these sub-symmetries, the solver needs to detect when it is in a subproblem where the sub-symmetry is active.
    Then, the sub-symmetry must be handled, which can be done using methods similar to those for handling global symmetries.

    As observed by~\cite{bendotti2021orbitopal}, when handling the sub-symmetries of all solution subsets in~$\mathbb S$ simultaneously, one needs to be slightly more careful as for global symmetries, however.
    The general idea remains the same: we want to disregard solutions such that for each orbit at least one representative remains.
    For a solution subset~$Q_s$, let~$\sigma^s_k$ for~$k \in [o_s]$ denote the orbits defined by $G_{Q_s}$.
    Furthermore, let~$O = \set{\sigma_k^s | s \in [q],\, k \in [o_s]}$ denote the family of all orbits of the considered solution subsets~$Q_s$.
    Notice that~$O$ does not necessarily partition the solutions of the solution subsets, as the orbits of different solution subsets may overlap.
    As a consequence, we need to be more careful in choosing the representative solution~$r(\sigma)$ for each orbit~$\sigma \in O$, as for a given orbit~$\sigma$ the set~$\sigma \setminus \set{r(\sigma)}$ may contain a representative of another orbit~$\sigma' \in O$.
    Bendotti et al.~\cite{bendotti2021orbitopal} propose to choose \emph{orbit-compatible representatives} that ensure that there always remains at least one solution in every orbit, when restricting to generalized representatives.

    If one only handles symmetries arising from permutations of columns in the matrix of binary variables, the following structure in the set of sub-symmetries~$\mathbb S$ ensures that generalized representatives can be defined easily~\cite{bendotti2021orbitopal}.
    Let~$Q \in \mathbb S$.
    For a solution matrix~$x \in Q$, let~$x(R, C)$ denote the submatrix of~$x$ obtained by restricting to rows~$R \subseteq [m]$ and columns~$C \subseteq [n]$.
    The symmetry group~$G_Q$ is the \emph{sub-symmetric group} with respect to~$(R, C)$ if it contains all the permutations of the columns of~$x(R, C)$.
    If~$G_Q$ is the sub-symmetric group, then~$Q$ is called \emph{sub-symmetric} with respect to $(R, C)$.
    Now, let~$\mathbb S$ be a set of solution subsets such that every~$Q_s \in \mathbb S$ is sub-symmetric with respect to~$(R_s, C_s)$.
    For every orbit~$\sigma_s^i$ of~$G_{Q_s}$, choose the representative~$x_s^i \in \sigma_s^i$ such that the submatrix~$x_s^i(R_s, C_S)$ is lexicographically maximal in its orbit, i.e., its columns are lexicographically non-increasing.
    Then, these representatives are orbit-compatible~\cite{bendotti2021orbitopal}.

    For the MKP example, the capacity sub-symmetries arise in the solution subsets~$\mathbb S_\text{MKP}$.
    A solution subset~$Q^i_{j,j'} \in \mathbb S_\text{MKP}$ is sub-symmetric with respect to~$(\set{i, \dots, m}, \set{j, j'})$.
    Whenever the sub-symmetry is active, one can thus handle it by enforcing that the columns of the submatrix~$y(\set{i, \dots, m}, \set{j, j'})$ are lexicographically non-increasing.

    \subsection{Handling Sub-symmetries}\label{sec:handlesubsym}
    An approach to handle a sub-symmetry is to add sub-symmetry-handling inequalities to the model~\cite{bendotti2020symmetry}.
    We briefly describe the framework of~\cite{bendotti2020symmetry}.
    Let~$Q_s$ be a solution subset that is sub-symmetric with respect to~$(R_s, C_s)$.
    We introduce an integer variable~$z_s$ such that~$z_s = 0$ if and only if~$x \in Q_s$ and $z_s \ge 1$ otherwise.
    Let~$c_{j},c_{j+1}$ be two consecutive columns in the submatrix~$x(R_s, C_s)$.
    Then, the \emph{partial sub-symmetry-handling inequality} is
    \begin{equation} \label{eq:partsubsym}
        x_{r_1,c_{j+1}} \le z_s + x_{r_1,c_j}\qquad \text{where $r_1 = \min R_s$.}
    \end{equation}
    If~$z_s = 0$, \eqref{eq:partsubsym} ensures that, for all pairs of consecutive columns, the first row of the submatrix~$x(R_s, C_s)$ is lexicographically non-increasing.
    Otherwise, if~$z_s \ge 1$, the inequality is trivially satisfied.

    In the MKP~example, the auxiliary variable~$z$ for~$Q^i_{j,j'} \in \mathbb S_\text{MKP}$ can be expressed as
    \begin{equation} \label{eq:actMKP}
        z = \left|c_j - c_{j'} - \sum_{k=1}^{i-1} w_k (y_{k,j} - y_{k,j'})\right|.
    \end{equation}
    To express~$z$ in the IP formulation, we need to linearize~\eqref{eq:actMKP}.
    Therefore, we write~$\alpha = c_j - c_{j'} - \sum_{k=1}^{i-1} w_k (y_{k,j} - y_{k,j'})$ for brevity, and introduce non-negative variables~$\alpha^+,\alpha^- \ge 0$, binary variables~$z^+,z^- \in \set{0,1}$, and constraints
    \begin{equation}
        \alpha^+ \le Mz^+,\, \alpha^- \le M z^-,\, \alpha^+ + \alpha^- \ge z^+ + z^-,\, \alpha = \alpha^+ - \alpha^-,\, z^+ + z^- \le 1.
    \end{equation}
    Here, $M$ is a sufficiently large constant, e.g.,~$\max\{c_j, c_{j'}\} + \sum_{k=1}^{i-1} w_k$.
    Then, $z = z^+ + z^-$ indicates whether the sub-symmetry is active.
    Note that the constraints indeed ensure that~$|\alpha| = \alpha^+ + \alpha^- = 0$ if and only if~$z = 0$.
    The sub-symmetry-handling inequality thus is
    \begin{equation}
        y_{i,j'} \le z + y_{i,j}.
    \end{equation}

    Notice that the Inequalities~\eqref{eq:partsubsym} are not sufficient to fully handle the symmetry.
    Indeed, they only ensure that the first row of the submatrix is lexicographically ordered.
    If the entries on the first row are equal, subsequent rows need to be considered until a \emph{tie-break} row is found.
    To this end, the set~$\mathbb S$ is extended to~$\tilde{\mathbb S}$ with additional \emph{tie-break subsets} for which Inequalities~\eqref{eq:partsubsym} break the symmetry on rows where the previous rows have equal entries.
    The set~$\tilde{\mathbb S}$ can be defined as
    \begin{equation}
        \tilde{\mathbb S} = \set{\tilde Q_s(i, j) | \text{$s \in [q]$, $i \in \set{1, \dots, {|R_s|}}$, $j \in \set{2, \dots, {|C_s|}}$}}
    \end{equation}
    where
    \begin{equation}
        \tilde Q_s(i, j) = \set{x \in Q_s | \text{$x_{r,c^s_{j-1}} = x_{r,c^s_j}$ for all $r \in \set{r^s_1, \dots, r^s_{i-1}}$}}.
    \end{equation}
    For our MKP example, one can verify that the set~$\mathbb S_\text{MKP}$ already includes the tie-break sets for every sub-symmetry.
    Therefore, it is not necessary to extend~$\mathbb S_\text{MKP}$ with additional tie-break subsets.

    \section{Activation Handler}
    \label{sec:ah}

    The existing method of handling sub-symmetries with inequalities has a number of limitations.
    For every sub-symmetry that we want to handle, it is necessary to add explicit SHIs to the formulation, leading to a blow-up of the IP\@.
    The size of the formulation increases even more with the addition of tie-break sets, and for problems where additional variables or constraints are necessary to express the auxiliary~$z$-variable in the formulation.
    Additionally, the inequalities are rather weak for symmetry handling.
    In particular, the variable-based approach is not immediately able to activate more sophisticated symmetry-handling methods such as orbitopal fixing.

    To circumvent these issues, we introduce a new approach for handling sub-symmetries.
    In our framework, we decouple the activation of sub-symmetries from the explicit IP~formulation itself, and instead adopt a more flexible approach.
    The modeler of the problem can provide a set of rules that define when a sub-symmetry becomes active to the IP~solver.
    We call this set of rules the \emph{activation handler}.
    The activation handler checks for a node~$a$ in the B\&B tree whether the rules hold.
    If that is the case, the activation handler activates a symmetry-handling method in the solver, to handle the sub-symmetry.

    In our framework, both the activation and handling of sub-symmetries is flexible.
    The modeler can describe the rules for activation with any custom implementation that uses information from the current state of the solver, such as variable fixings at the current node of the B\&B tree, to determine whether a subproblem with sub-symmetry is active.
    Furthermore, the activation handler can be used with any symmetry-handling method from the literature, and is not restricted to inequality-based approaches.
    Neither activation nor symmetry handling needs to be encoded in the formulation directly, keeping the IP compact.

    As the rules-based approach of activation handlers is rather generic, we will illustrate the concept for several problems in Sec.~\ref{sec:applications}.
    First, we describe sub-symmetries that arise in certain applications and compare how the SHI-based approach and activation handlers can be used to handle them.
    Afterwards, we compare the numerical performance of the different approaches in Sec.~\ref{sec:experiments}.
    Our implementation of the activation handler framework is publicly available\footnote{See \url{https://github.com/stenwessel/activation-handler}.}, as well as the setup of our experiments.
    The generic framework can easily be adapted by practitioners for use in other applications.

    \section{Application}
    \label{sec:applications}
    In this section, we discuss how sub-symmetries arise for three types of problems: the multiple knapsack problem, the unit commitment problem, and the maximum $k$-colorable subgraph problem.
    We describe how SHIs can be applied to handle sub-symmetries, as well as the activation handler framework.
    The activation handler uses information from the solver about variable fixings at a node of the B\&B tree.
    To this end, we define for a node~$a$ of the B\&B tree the sets~$F^a_0$ and~$F^b_1$, which denote the variables that are fixed to~$0$ or~$1$ at node~$a$, respectively.

    \subsection{Multiple Knapsack Problem}
    We introduced the MKP, its symmetries, and the capacity sub-symmetries in Example~\ref{ex:mkp}.
    Notice that the number of solution subsets in~$\mathbb S_\text{MKP}$ is rather large, as we consider all pairs of knapsacks.
    When handling the symmetry with inequalities, considering all solution subsets is intractable, as potentially every subset of knapsacks might define a sub-symmetry.
    That is, exponentially many SHIs as well as auxiliary variables and constraints need to be added to the problem.
    Therefore, we only consider SHIs for \emph{consecutive} pairs of knapsacks, i.e.,~$j' = j + 1$.
    We hence add~$\bigO(mn)$~SHIs to the formulation, with for each SHI four auxiliary variables and five auxiliary constraints.

    \paragraph{Sub-symmetry handling with activation handler}
    When handling sub-symmetries via activation handlers, we are more flexible in implementing the activation rules.
    Instead of enumerating every solution subset~$Q^i_{j_1,j_2}$ separately and checking if the sub-symmetry is active, we can use a single activation handler in the model.
    The activation handler returns all submatrices of~$y$ that contain active sub-symmetries at a given node of the B\&B tree.

    The activation handler identifies whether the placement of items~$\set{1, \dots, i-1}$ is fixed at node~$a$, according to the variable fixings~$F^a_0$ and~$F^a_1$.
    For every item~$i$ for which the previous holds, the activation handler checks whether there are knapsacks of equal remaining capacity, after placement of items~$\set{1, \dots, i-1}$.
    Suppose that for item~$i$ the knapsacks~$j_{k_1}, \dots, j_{k_r}$ have equal remaining capacity.
    Then, the activation handler reports the submatrix~$y(\set{i, \dots, m}, \set{j_{k_1}, \dots, j_{k_r}})$, for which the capacity sub-symmetry is now active at node~$a$.
    In this way, finding all active capacity sub-symmetries is linear in the number of variables in the matrix~$y$, as we can simply perform a linear scan over the rows of~$y$ and checking the variable fixings.
    Note that checking whether a variable is fixed can be done in constant time, as this information is available from the solver.

    The activated submatrices are then passed to a high-level symmetry-handling constraint in the solver.
    Several methods can be used to handle symmetry in the submatrix.
    In our implementation, we use orbitopal fixing for packing orbitopes~\cite{KaibelEtAl2011} to handle the active sub-symmetries.
    However, when all knapsacks are in the activated submatrix (i.e., all columns), the sub-symmetry is instead handled with orbitopal fixing for the stronger partitioning orbitope.

    \subsection{Unit Commitment Problem}
    Another problem in which we can handle sub-symmetries is the min-up/min-down unit commitment problem~(MUCP), as introduced in~\cite{bendotti2020symmetry}.
    We are given a set of production units~$\mathcal U$ with~$|\mathcal U| = n$, and a discrete time horizon~$\mathcal T = \set{1, \dots, T}$ for which a certain non-negative demand~$D_t$ needs to be satisfied at every time~$t \in \mathcal T$.
    Every production unit~$j \in \mathcal U$ can be either \emph{up} or \emph{down} at every~$t \in \mathcal T$.
    When a unit is up, its production is between a minimum and maximum production capacity~$P_\text{min}^j$ and~$P_\text{max}^j$, and it must remain up for at least~$L^j$~time steps.
    When a unit is down, its production is zero and it must remain down for at least~$\ell^j$~time steps.
    We furthermore have for every unit~$j$ a start-up cost~$c_0^j$, a fixed cost~$c_f^j$ for every time step the unit is up, and a production cost~$c_p^j$ proportional to its production.
    The goal is to find a \emph{production schedule} satisfying the production demand at every time step and the min-up and min-down constraints, while minimizing the total cost.

    Let the variables~$x_{t,j} \in \set{0, 1}$ indicate whether unit~$j \in \mathcal U$ is up at time~$t \in \mathcal T$,
    and~$u_{t,j} \in \set{0,1}$ whether unit~$j$ starts up at time~$t$.
    We omit further details of the IP~formulation we use for this problem, as they are not relevant for symmetry handling, and refer to~\cite{bendotti2020symmetry} for details.
    Let~$\mathcal X_\text{MUCP}$ denote the set of matrices~$(x_{t,j})$ that are feasible.
    Notice that the solution matrix~$x$ completely characterizes a solution, as the corresponding matrix~$u$ can be derived completely from~$x$.

    \paragraph{Symmetry in the MUCP}
    Symmetries are present globally in the MUCP when production units have identical properties, i.e., units where all of the properties~$(P_\text{min}, P_\text{max}, L, \ell, c_0, c_f, c_p)$ are equal.
    To make this more explicit, we partition the production units into $H$~\emph{types}, where a type~$h \in \set{1, \dots, H}$ consists of~$n_h$ identical units that we denote by~$\mathcal U_h = \set{j^h_1, \dots, j^h_{n_h}}$.
    For a type~$h$, we slightly abuse notation to denote its properties as $(P_\text{min}^h, P_\text{max}^h, L^h, \ell^h, c_0^h, c_f^h, c_p^h)$.
    We can then also partition the matrix variable~$x$ into~$H$~matrices~$x^h = (x_{t,j})_{t \in \mathcal T, j \in \mathcal U_h}$ for every type of production unit.
    The production units within each type are identical, and we can hence permute their production schedules.
    This corresponds to permuting the columns of~$x^h$.
    One possible way of breaking the symmetry is to restrict $x^h$ to the full orbitope for binary matrices of size~$T \times n_h$, i.e., by imposing that the columns of~$x^h$ are lexicographically non-increasing.

    The MUCP also exhibits sub-symmetries, as introduced in~\cite{bendotti2020symmetry}.
    Call a production unit~$j \in \mathcal U$ \emph{ready to start up} at some time~$t \in \mathcal T$ if the unit has been down continuously for at least the minimum downtime~$\ell^j$.
    In other words, when~$x_{t',j} = 0$ for all~$t' = t-\ell^j, \dots, t-1$ and $t \ge \ell^j + 1$.
    Now, suppose there are at least two units~$j_1, \dots, j_k \in \mathcal U^h$ of type~$h$ that are all ready to start up at some time~$t$.
    Then, their production schedules can be permuted from time~$t$ onwards, regardless of their schedule up to time~$t$.
    This thus defines a sub-symmetry where the columns of the submatrix~$x(\set{t, \dots, T}, \set{j_1, \dots, j_k})$ can be permuted.
    Analogously, one can identify sub-symmetries for two units ready to shut down at some time~$t \in \mathcal T$.
    These sub-symmetries are referred to as the \emph{start-up} and \emph{shut-down sub-symmetries}, respectively.

    \paragraph{Sub-symmetry-handling inequalities}
    Following the approach in~\cite{bendotti2020symmetry}, the start-up sub-symmetries can be handled with inequalities as follows.
    The handling of shut-down sub-symmetries is analogous, and we omit the details here.
    Let~$j_k^h, j^h_{k+1}$ be a pair of consecutive units of the same type~$h$.
    For brevity, let~$j = j_k^h$ and~$j' = j^h_{k+1}$.
    Then, the solution subsets
    \begin{equation}
        \check Q_{k,h}^t = \set{x \in \mathcal X_\text{MUCP} | \text{$x_{t',j} = x_{t',j'} = 0$ for all $t' = t - \ell^h, \dots, t-1$}}
    \end{equation}
    for all~$t \ge \ell^h + 1$, define when the start-up sub-symmetries occur.
    Note that no additional tie-break sets are necessary.
    For~$\check Q_{k,h}^t$, the corresponding auxiliary variable can be expressed as~$z = \sum_{t'=t-\ell^h}^{t-1} [x_{t',j} + x_{t',j'}]$, leading to the SHI
    \begin{equation} \label{eq:shi:mucp}
        x_{t,j'} \le z + x_{t,j}.
    \end{equation}
    Note that the~$z$-variable has a linear description in~$x$, and hence it is not necessary to add~$z$ as a new variable to the formulation.
    Instead, we can simply replace~$z$ directly with its linear expression in the SHI~\eqref{eq:shi:mucp}.

    The SHIs for the start-up sub-symmetries can be slightly strengthened to
    \begin{equation} \label{eq:strengthenedSHI}
        u_{t,j'} \le x_{t-\ell^h,j} + x_{t,j} + \sum_{t'=t-\ell^h + 1}^{t-1} u_{t',j},
    \end{equation}
    leading to a stronger LP~relaxation.
    A similar inequality can be obtained for the shut-down sub-symmetries,
    see~\cite{bendotti2020symmetry} for the derivation of the strengthened SHIs.

    \paragraph{Sub-symmetry-handling with activation handler}
    Handling sub-symmetry with an activation handler is similar to the approach for the MKP\@.
    We add a single activation handler to the model, that identifies all submatrices corresponding to active sub-symmetries in the following manner.
    For sake of presentation, we assume that all production units~$\mathcal U$ have the same type.
    In the more general case where we have multiple types of production units, we can simply apply our method to the unit types separately.

    Let~$a$ be a node of the B\&B tree.
    Define for every~$t \in \set{\ell + 1, \dots, T}$,
    \begin{equation}
        \check S^a_t = \set{j \in \mathcal U | \text{$x_{t',j} \in F^a_0$ for all $t' \in \set{t - \ell, \dots, t-1}$}}.
    \end{equation}
    That is, $\check S^a_t$ are the production units that are \emph{fixed} to be ready to start up at time~$t$ at node~$a$.
    For every subset~$\check S^a_{t}$ for which $|S^a_{t}| \ge 2$, the corresponding start-up sub-symmetry becomes active.
    Hence, the symmetry corresponds to column permutations of the submatrix~$x(\set{t, \dots, T},\check S^a_{t})$.
    We then use orbitopal fixing for full orbitopes to handle the sub-symmetry in the activated submatrix.

    Notice that we can find the units that are ready to start up for every time~$t \in \mathcal T$ in~$\bigO(nT)$~time, by iterating over the time horizon and a dynamic-programming approach.
    The shut-down sub-symmetries are activated and handled with an analogous approach.

    \subsection{Graph Coloring Problem}
    In this section, we discuss another application of the sub-symmetry-handling framework on a variant of the graph coloring problem.
    Consider an undirected graph~$G = (V, E)$ with~$|V| = n$ and let~$k$ be a positive integer.
    A \emph{$k$-coloring} of the graph is a function~$c\colon V \to [k]$ that assigns colors to the vertices, such that for any pair~$\set{i,j} \in E$ it holds that~$c(i) \neq c(j)$.
    In the \emph{max-$k$-colorable subgraph problem}~(MKCS), we want to find a subset~$V' \subseteq V$ of vertices of maximum size, such that~$V'$ induces a subgraph that admits a~$k$-coloring.
    Let the binary variables~$x_{i,r} \in \set{0,1}$ indicate that vertex~$i \in V$ is colored with color~$r \in [k]$.
    We then have the standard IP~formulation~\cite{JanuschowskiPfetsch2011}
    \begin{subequations} \label{eq:ip:mkcs}
        \begin{alignat}{5}
            \text{maximize}\quad && \sum_{i\in V}\sum_{r=1}^k x_{i,r} \span \label{eq:ip:mkcs:obj} \\
            \text{subject to}\quad && x_{i,r} + x_{j,r} &\le 1 &&\qquad \forall_{\set{i,j} \in E},\ \forall_{r \in [k]}, \label{eq:ip:mkcs:edge} \\
            && \sum_{r=1}^k x_{i,r} &\le 1 &&\qquad\forall_{i \in V}, \label{eq:ip:mkcs:vertex}\\
            && x_{i,r} &\in \set{0, 1} &&\qquad\forall_{i \in V},\ \forall_{r \in [k]}. \label{eq:ip:mkcs:x}
        \end{alignat}
    \end{subequations}
    Let~$\X_\text{MKCS}$ denote the set of all feasible solution matrices~$x$.
    The MKCS~problem is NP-hard for any fixed~$k$, with~$1 \le k < n$~\cite{JanuschowskiPfetsch2011}.

    The type of symmetry we consider in this problem is the equivalence between the colors.
    Indeed, in any feasible solution, the color indices can be permuted to obtain a different, equivalent, feasible solution.
    This corresponds to the permutation of the columns in the solution matrix~$(x_{i,r})_{i\in V,r \in [k]}$.

    \paragraph{Sub-symmetry in the MKCS problem}
    The formulation~\eqref{eq:ip:mkcs} also exhibits sub-symmetries, c.f.~\cite{bendotti2020symmetry}.
    Consider two distinct colors~$c_1,c_2 \in [k]$ and a subset of the vertices~$R \subseteq V$ such that the neighbors of~$R$, denoted by~$N(R)$, are colored with neither~$c_1$ nor~$c_2$.
    Then, the colors~$c_1$ and~$c_2$ can be permuted within~$R$.
    Hence, the sub-symmetry occurs within the solution subset
    \begin{equation}
        Q^R_{c_1,c_2} = \set{x \in \X_\text{MKCS} | \text{$x_{i,c_1} = x_{i,c_2} = 0$ for all $i \in N(R)$}}
    \end{equation}
    with the sub-symmetry corresponding to the permutations of the columns in the submatrix~$x(R,\set{c_1,c_2})$.
    Notice that there are an exponential number of such solution subsets, because of the exponential number of subsets of vertices.
    Handling all such sub-symmetries with inequalities is therefore infeasible in practice, because of the blow-up in the number of auxiliary constraints.
    The adopted approach in~\cite{bendotti2020symmetry} is to only consider some subsets of vertices to handle sub-symmetries, leading to a quadratic number of SHIs.
    Within our activation handler framework, we can instead approach this differently and dynamically find relevant vertex subsets at the nodes of the B\&B tree.

    \paragraph{Sub-symmetry-handling with activation handler}
    Let~$a$ be a node of the B\&B tree and fix a pair of distinct colors~$c_1,c_2$.
    Let~$S \subseteq V$ be a set of vertices, such that~$i \in S$ when~$x_{i,c_1},x_{i,c_2} \in F^0_a$.
    Let now~$R \define V \setminus S$, and notice that indeed~$N(R)$ only contains vertices that are fixed to be not colored by~$c_1$ nor~$c_2$.
    In particular, the connected components~$R_1, \dots, R_\ell$ of the subgraph induced by~$R$ all satisfy this property.
    Hence, we can handle sub-symmetry in the submatrices~$x(R_j,\set{c_1,c_2})$ for all connected components~$R_j$.

    An activation handler can implement this activation rule via
    an iterated depth-first-search approach that efficiently identifies the connected components.
    For non-trivial connected components, the sub-symmetry is handled by orbitopal fixing for packing orbitopes on the identified submatrix.
    In particular, in contrast to the SHI-based approach, the activation handler can handle \emph{all} possible color-pair-based sub-symmetries without increasing the size of the IP formulation.

    \section{Experimental Results}\label{sec:experiments}
    In this section, we compare the sub-symmetry-handling methods using experiments on instances of the problems introduced above.

    \subsection{Instances}
    For the MKP, we generate random instances in the four standard classes of problems from the literature on MKP~\cite{MARTELLO1987213,PISINGER1999528,fukunaga2011branch}:
    \begin{itemize}
        \item \emph{uncorrelated}, where the weights~$w_i$ and profits~$p_i$ are uniformly distributed integers in the closed range~$[\ell, L]$,
        \item \emph{weakly correlated}, where the weights~$w_i$ are uniformly distributed integers in the closed range~$[\ell, L]$ and~$p_i$ are uniformly distributed integers in the closed range~$[\max\set{1, w_i - (L - \ell)/10}, w_i + (L - \ell)/10]$,
        \item \emph{strongly correlated}, where the weights~$w_i$ are uniformly distributed integers in the closed range~$[\ell, L]$ and~$p_i = w_i + (L - \ell)/10$,
        \item \emph{multiple subset-sum}, where the weights~$w_i$ are uniformly distributed integers in the closed range~$[\ell, L]$ and~$p_i = w_i$.
    \end{itemize}
    We generate our instances with~$\ell = 10$, $L = 1000$.
    The capacity of every knapsack is set to~$c_j = \lfloor \frac 1 2 \sum_{i=1}^m w_i/n \rfloor$, such that the total capacity is approximately half of the total weight of all items.
    To introduce symmetry in the problem, we generate multiple items with the same weight with an approach similar to Bendotti et al.~\cite{bendotti2020symmetry}.
    Generated weights are duplicated~$d$ times, where~$d$ is a uniformly random integer in~$[1, fm]$, where we call~$f \in \set{\frac 1 2, \frac 1 3, \frac 1 4, \frac 1 8}$ the \emph{symmetry factor}.
    Larger values of~$f$ generate larger groups of items with equal weight, leading to a more symmetric instance.
    For generating the profit values for the items within an equal-weight group, we consider two types of instances:
    \begin{itemize}
        \item \emph{equal profit}, where every item in the equal-weight group also has equal profit, generated according to the item class above,
        \item \emph{free profit}, where every item in the equal-weight group has a profit value generated according to the item class above.
    \end{itemize}
    Note that for the \emph{strongly correlated} and \emph{multiple subset-sum} classes, we only generate instances for \emph{equal profit} as both types of instances are equivalent.
    We generate groups of equal-weight items until we have generated~$m$ items.
    For every pair of~$(m, n) \in \set{(48, 12), (60, 10), (60, 30), (75, 15), (100, 10)}$, we generate~$20$ instances for every combination of item class, symmetry factor, and duplication type, yielding a total of $2400$~instances.

    For the MUCP, we use the same generated instances that are used for the experimental evaluation in~\cite{bendotti2020symmetry}, of which we are grateful to the authors for providing them.
    For the MKCS problem, we use the standard DIMACS Color02 set of graph coloring instances\footnote{Obtained from~\url{https://mat.tepper.cmu.edu/COLOR02/}.}, with~$k \in \set{5,6,8,10}$.

    \subsection{Experimental setup}
    All experiments are run with the development version of SCIP~7.0.3 (Git hash \texttt{3671128c}) with the SoPlex LP solver (Release~600)~\cite{scip}, on a single core of an Intel Xeon Platinum~8260~CPU running at~$2.4$~GHz, with a memory limit set at $10$~GB of RAM\@ and a solving-time limit of~$3600$~seconds for MKP and MUCP, and $7200$~seconds for the MKCS problem.
    The IP model is constructed in Python~3.10 using the PySCIPOpt interface that exposes the SCIP API in Python.
    The activation handler is implemented in SCIP as a new plugin, and can be added to the model with the PySCIPOpt interface.

    Every instance is solved with five different settings, in order to compare performance of the different symmetry-handling methods:
    \begin{itemize}
        \item No-Sym: Formulation with SCIP internal symmetry handling turned off.
        \item Default: Formulation with SCIP default parameters.
        \item Orbitope: Formulation with orbitope constraints for (global) symmetry handling.
        \item Ineq: Formulation with SHIs.
        \item Act: Formulation with orbitope constraints for (global) symmetry handling and activation handler for sub-symmetries.
    \end{itemize}
    For the MKP, all models except for No-Sym include orbitope constraints for handling symmetry between identical items.
    In the orbitopes for symmetries between identical items and symmetries between identical knapsacks, we use a compatible ordering of the variables such that orbitopal fixing for all orbitopes can be performed simultanuously, without introducing any conflicts.
    For the MUCP, we use the strengthened SHIs~\eqref{eq:strengthenedSHI} in the Ineq model.
    For MKCS, we consider two variants of the Act model: Act-AllPairs and Act-Consec which activate the sub-symmetries for all pairs of colors, or only pairs of consecutive colors, respectively.
    We do not consider a formulation with SHIs for MKCS, to avoid making arbitrary choices for which vertex subsets to consider.
    The orbitope constraints in SCIP use orbitopal fixing, as discussed in Sec.~\ref{sec:introduction}.

    \subsection{Results}
    The results are summarized in Tables~\ref{tab:mkp}, \ref{tab:mucp}, \ref{tab:mkcs:small}, and \ref{tab:mkcs:large} for the MKP, MUCP, MKCS with~$k=5,6$, and MKCS with~$k=8,10$, respectively.
    In the reported results, we summarize the results for instances in classes, based on the solving time of the tested models.
    We use the notation~$[a, b)$ to denote the set of instances for which all models have a solving time of at least~$a$ and below~$b$~seconds.
    We exclude instances from our test set where all models reach the time limit.
    For every instance class, we report the number of instances in the class~(\#).
    For every model, the number of instances solved to optimality (Opt) is reported, as well as the mean solving time, in seconds, of all instances.
    The mean solving time is the shifted geometric mean, with a shift of~$1$~second.
    For instances that are not solved to optimality within the time limit, the solving time is set to~$3600$~seconds for MKP and MUCP, and $7200$~seconds for the MKCS~problem.

    \begin{table}[p]
        \centering
        \caption{Summarized numerical results for MKP\@.
        For $532$~instances, all models reach the time limit and are excluded from the table.} \label{tab:mkp}
        \scriptsize
        \begin{tabular}{l *{11}{@{\hskip 1em}r}}
            \toprule
             && \multicolumn{2}{@{\hskip .5em}c@{\hskip 1em}}{No-Sym} & \multicolumn{2}{@{\hskip .5em}c@{\hskip 1em}}{Default} & \multicolumn{2}{@{\hskip .5em}c@{\hskip 1em}}{Ineq} & \multicolumn{2}{@{\hskip .5em}c@{\hskip 1em}}{Orbitope} & \multicolumn{2}{@{\hskip .5em}c@{\hskip 1em}}{Act} \\
            \cmidrule(r){3-4} \cmidrule(r){5-6} \cmidrule(r){7-8} \cmidrule(r){9-10} \cmidrule{11-12}
            Instances       & {$\#$} & {Opt} & {Time} & {Opt} & {Time} & {Opt} & {Time} & {Opt} & {Time} & {Opt} & {Time} \\
            \midrule
            All             & 1868   & 1166  & 52.7   & 1387  & 35.3   & 625   & 1122.9 & 1636  & 17.6   & 1819  & 13.9   \\
            $[0,100)$       & 277    & 277   & 0.0    & 277   & 0.2    & 277   & 14.4   & 277   & 0.1    & 277   & 0.2    \\
            $[100,1800)$    & 269    & 269   & 0.4    & 269   & 0.5    & 269   & 370.5  & 269   & 0.3    & 269   & 0.4    \\
            $[1800,\infty)$ & 1322   & 620   & 256.8  & 841   & 139.1  & 79    & 3460.2 & 1090  & 56.0   & 1273  & 40.3   \\
            \bottomrule
        \end{tabular}
    \end{table}
    \begin{table}[p]
        \centering
        \caption{Summarized numerical results for MUCP\@.
        For $12$~instances, all models reach the time limit and are excluded from the table.} \label{tab:mucp}
        \scriptsize
        \begin{tabular}{l *{9}{@{\hskip 1em}r}}
            \toprule
             && \multicolumn{2}{@{\hskip .5em}c@{\hskip 1em}}{No-Sym} & \multicolumn{2}{@{\hskip .5em}c@{\hskip 1em}}{Default} & \multicolumn{2}{@{\hskip .5em}c@{\hskip 1em}}{Ineq} & \multicolumn{2}{@{\hskip .5em}c@{\hskip 1em}}{Act} \\
            \cmidrule(r){3-4} \cmidrule(r){5-6} \cmidrule(r){7-8} \cmidrule{9-10}
            Instances      & {$\#$} & {Opt} & {Time} & {Opt} & {Time} & {Opt} & {Time} & {Opt} & {Time} \\
            \midrule
            All            & 268    & 172   & 162.7  & 216   & 70.8   & 240   & 131.4  & 266   & 39.0   \\
            $[0,10)$       & 26     & 26    & 3.1    & 26    & 3.0    & 26    & 5.0    & 26    & 3.6    \\
            $[10,300)$     & 101    & 101   & 19.2   & 101   & 14.1   & 101   & 31.3   & 101   & 15.8   \\
            $[300,\infty)$ & 34     & 1     & 3353.4 & 7     & 2789.0 & 18    & 1972.0 & 32    & 751.3  \\
            \bottomrule
        \end{tabular}
    \end{table}
    \begin{table}[p]
        \centering
        \caption{Summarized numerical results for MKCS instances with~$k = 5, 6$.
        For $92$~instances, all models reach the time limit and are excluded from the table.} \label{tab:mkcs:small}
        \scriptsize
        \begin{tabular}{l *{11}{@{\hskip 1em}r}}
            \toprule
            && \multicolumn{2}{@{\hskip .5em}c@{\hskip 1em}}{No-Sym} & \multicolumn{2}{@{\hskip .5em}c@{\hskip 1em}}{Default} & \multicolumn{2}{@{\hskip .5em}c@{\hskip 1em}}{Orbitope} & \multicolumn{2}{@{\hskip .5em}c@{\hskip 1em}}{Act-Consec} & \multicolumn{2}{@{\hskip .5em}c@{\hskip 1em}}{Act-AllPairs} \\
            \cmidrule(r){3-4} \cmidrule(r){5-6} \cmidrule(r){7-8} \cmidrule(r){9-10} \cmidrule{11-12}
            Instances      & {$\#$} & {Opt} & {Time} & {Opt} & {Time} & {Opt} & {Time} & {Opt} & {Time} & {Opt} & {Time} \\
            \midrule
            All            & 147    & 139   & 18.4   & 145   & 17.4   & 141   & 20.4   & 141   & 22.0   & 141   & 21.9   \\
            $[0,10)$       & 55     & 55    & 0.4    & 55    & 0.5    & 55    & 0.7    & 55    & 0.7    & 55    & 0.9    \\
            $[10,600)$     & 63     & 63    & 26.1   & 63    & 24.9   & 63    & 31.8   & 63    & 36.2   & 63    & 35.0   \\
            $[600,\infty)$ & 13     & 7     & 3482.0 & 11    & 3263.3 & 8     & 2990.4 & 8     & 3660.4 & 8     & 3703.8 \\
            \bottomrule
        \end{tabular}
    \end{table}
    \begin{table}[p]
        \centering
        \caption{Summarized numerical results for MKCS instances with~$k = 8, 10$.
        For $91$~instances, all models reach the time limit and are excluded from the table.} \label{tab:mkcs:large}
        \scriptsize
        \begin{tabular}{l *{11}{@{\hskip 1em}r}}
            \toprule
             && \multicolumn{2}{@{\hskip .5em}c@{\hskip 1em}}{No-Sym} & \multicolumn{2}{@{\hskip .5em}c@{\hskip 1em}}{Default} & \multicolumn{2}{@{\hskip .5em}c@{\hskip 1em}}{Orbitope} & \multicolumn{2}{@{\hskip .5em}c@{\hskip 1em}}{Act-Consec} & \multicolumn{2}{@{\hskip .5em}c@{\hskip 1em}}{Act-AllPairs} \\
            \cmidrule(r){3-4} \cmidrule(r){5-6} \cmidrule(r){7-8} \cmidrule(r){9-10} \cmidrule{11-12}
            Instances      & {$\#$} & {Opt} & {Time} & {Opt} & {Time} & {Opt} & {Time} & {Opt} & {Time} & {Opt} & {Time} \\
            \midrule
            All            & 147    & 140   & 7.6    & 143   & 9.0    & 138   & 16.9   & 136   & 20.1   & 139   & 25.2   \\
            $[0,10)$       & 60     & 60    & 0.1    & 60    & 0.3    & 60    & 0.5    & 60    & 0.8    & 60    & 1.1    \\
            $[10,600)$     & 60     & 60    & 7.9    & 60    & 12.1   & 60    & 26.8   & 60    & 36.2   & 60    & 45.7   \\
            $[600,\infty)$ & 7      & 5     & 4810.1 & 3     & 5452.5 & 4     & 6204.7 & 3     & 5284.1 & 4     & 4515.3 \\
            \bottomrule
        \end{tabular}
    \end{table}

    In Table~\ref{tab:mkp}, we see that over all instances, the activation handler method solves more instances to optimality within the time limit, compared to the other models.
    We can see that symmetry handling is highly relevant for the MKP problem.
    SCIP's state-of-the-art symmetry-handling methods reduce the running time by roughly~\SI{33}{\percent}.
    Our activation handler approach reduces the running time of this already very competitive setting by further~\SI{60}{\percent}.
    Comparing the different sub-symmetry-handling methods, the additional overhead necessary for the inequalities method is too large for handling this type of sub-symmetries.
    From the small and medium classes, we see there is a substantial difference between solving time for the inequalities method.
    There is also a slight improvement in running time for the global orbitope and activation handler methods, compared to default SCIP and the model where no symmetry handling is performed.
    For the large instances we see that the activation handler method solves more instances to optimality with a clear improvement in solving time.

    In Table~\ref{tab:mucp} we see similar results.
    Overall, SCIP's symmetry-handling methods improve the running time by roughly~\SI{55}{\percent}, whereas the activation handler reduces it even further by~\SI{45}{\percent}.
    We omit the Orbitope model in the results, as SCIP automatically finds these orbitopes in the Default model.
    The inequalities method, in contrast with the results for the MKP, shows improvement on the default models for the large instances, confirming the results of~\cite{bendotti2020symmetry}.
    This difference compared to the results for the MKP is likely caused by the auxiliary $z$-variables, that can here be expressed linearly with no additional constraints.
    The activation handler outperforms all other models;
    it solves considerably more large instances to optimality.

    For the MKCS instances in Table~\ref{tab:mkcs:small}, (sub-)symmetry handling on the colors seems to have an adverse effect on the performance, and the default SCIP model performs best on average.
    This is likely caused by SCIP's automatic symmetry detection, which might discover symmetries in the graph that are not handled in the global orbitope and activation handler models.
    In this case, it seems that graph symmetries lead to more powerful performance improvements for this problem.
    In Table~\ref{tab:mkcs:large}, the results are similar, but for the large instances the average solving times are increased for every model, compared to the instances in Table~\ref{tab:mkcs:small} where~$k$ is smaller.
    Here, we see that the activation handler outperforms only handling global symmetries in the orbitope model;
    although the model with no symmetry handling generally outperforms both methods.
    This would imply that for very large instances, the activation handler likely scales better with the instance size than only handling global symmetries.

    \section{Conclusion} \label{sec:conclusion}
    We have introduced a new framework for sub-symmetry handling in integer programming.
    The approach is flexible, as it allows modelers to generically implement rules that define sub-symmetries, and it can be used in conjunction with any kind of symmetry-handling method in the existing literature.
    We have shown via a number of applications that it can be used in various settings, and computational experiments show that our framework substantially improves on the performance of the state-of-the-art methods.

    For future work, we are interested in looking into the interface of the activation handler framework.
    Currently, its flexibility allows users to write custom code for the implementation of the activation rules.
    For some problems, the activation rules may have a common structure that can be exploited,
    e.g., the capacity sub-symmetries for the MKP can also be applied for other bin-packing or scheduling problems.
    Such activation handlers may be generalized to be applicable for a broader set of formulations, for which a \emph{domain specific language} might allow practitioners to more easily specify sub-symmetries in their models.

    \subsubsection{Acknowledgments} \label{subsubsec:acknowledgements}
    We thank the authors of~\cite{bendotti2020symmetry} for providing us with the MUCP~instances originally used in their experiments.

    \bibliographystyle{splncs04}
    \bibliography{references}

\begin{thebibliography}{10}
\providecommand{\url}[1]{\texttt{#1}}
\providecommand{\urlprefix}{URL }
\providecommand{\doi}[1]{https://doi.org/#1}

\bibitem{babai1983canonical}
Babai, L., Luks, E.M.: Canonical labeling of graphs. In: Proceedings of the
  fifteenth annual {ACM} symposium on Theory of computing - {STOC}
  {\textquotesingle}83. {ACM} Press (1983). \doi{10.1145/800061.808746}

\bibitem{bendotti2020symmetry}
Bendotti, P., Fouilhoux, P., Rottner, C.: Symmetry-breaking inequalities for
  {ILP} with structured sub-symmetry. Mathematical Programming
  \textbf{183}(1),  61--103 (2020). \doi{10.1007/s10107-020-01491-4}

\bibitem{bendotti2021orbitopal}
Bendotti, P., Fouilhoux, P., Rottner, C.: Orbitopal fixing for the full
  (sub-)orbitope and application to the unit commitment problem. Mathematical
  Programming  \textbf{186}(1),  337--372 (2021).
  \doi{10.1007/s10107-019-01457-1}

\bibitem{DoornmalenHojny2022}
van Doornmalen, J., Hojny, C.: Efficient propagation techniques for handling
  cyclic symmetries in binary programs (2022), available at
  \url{https://optimization-online.org/2022/03/8812/}

\bibitem{FischettiLiberti2012}
Fischetti, M., Liberti, L.: Orbital shrinking. In: Mahjoub, A.R., Markakis, V.,
  Milis, I., Paschos, V.T. (eds.) Combinatorial Optimization, Lecture Notes in
  Computer Science, vol.~7422, pp. 48--58. Springer Berlin Heidelberg (2012)

\bibitem{Friedman2007}
Friedman, E.J.: Fundamental domains for integer programs with symmetries. In:
  Dress, A., Xu, Y., Zhu, B. (eds.) Combinatorial Optimization and
  Applications, Lecture Notes in Computer Science, vol.~4616, pp. 146--153.
  Springer Berlin Heidelberg (2007). \doi{10.1007/978-3-540-73556-4_17},
  \url{http://dx.doi.org/10.1007/978-3-540-73556-4_17}

\bibitem{fukunaga2011branch}
Fukunaga, A.S.: A branch-and-bound algorithm for hard multiple knapsack
  problems. Annals of Operations Research  \textbf{184}(1),  97--119 (2011).
  \doi{10.1007/s10479-009-0660-y}

\bibitem{scip}
Gamrath, G., Anderson, D., Bestuzheva, K., Chen, W.K., Eifler, L., Gasse, M.,
  Gemander, P., Gleixner, A., Gottwald, L., Halbig, K., Hendel, G., Hojny, C.,
  Koch, T., Le~Bodic, P., Maher, S.J., Matter, F., Miltenberger, M.,
  M{\"u}hmer, E., M{\"u}ller, B., Pfetsch, M.E., Schl{\"o}sser, F., Serrano,
  F., Shinano, Y., Tawfik, C., Vigerske, S., Wegscheider, F., Weninger, D.,
  Witzig, J.: {The SCIP Optimization Suite 7.0}. ZIB-Report 20-10, Zuse
  Institute Berlin (2020),
  \url{http://nbn-resolving.de/urn:nbn:de:0297-zib-78023}

\bibitem{GareyJohnson1979}
Garey, M.R., Johnson, D.S.: Computers and Intractability: A Guide to the Theory
  of {NP}-Completeness. W. H. Freeman \& Co. (1979)

\bibitem{Hojny2020}
Hojny, C.: Packing, partitioning, and covering symresacks. Discrete Applied
  Mathematics  \textbf{283},  689--717 (2020).
  \doi{https://doi.org/10.1016/j.dam.2020.03.002},
  \url{http://www.sciencedirect.com/science/article/pii/S0166218X20300949}

\bibitem{hojny2019polytopes}
Hojny, C., Pfetsch, M.E.: Polytopes associated with symmetry handling.
  Mathematical Programming  \textbf{175}(1),  197--240 (2019).
  \doi{10.1007/s10107-018-1239-7}

\bibitem{JanuschowskiPfetsch2011}
Januschowski, T., Pfetsch, M.E.: The maximum $k$-colorable subgraph problem and
  orbitopes. Discrete Optimization  \textbf{8}(3),  478--494 (2011).
  \doi{http://dx.doi.org/10.1016/j.disopt.2011.04.002}

\bibitem{KaibelEtAl2011}
Kaibel, V., Peinhardt, M., Pfetsch, M.E.: Orbitopal fixing. Discrete
  Optimization  \textbf{8}(4),  595--610 (2011).
  \doi{http://dx.doi.org/10.1016/j.disopt.2011.07.001},
  \url{http://www.sciencedirect.com/science/article/pii/S1572528611000430}

\bibitem{kaibel2008packing}
Kaibel, V., Pfetsch, M.E.: Packing and partitioning orbitopes. Mathematical
  Programming  \textbf{114}(1),  1--36 (2008). \doi{10.1007/s10107-006-0081-5}

\bibitem{LandDoig1960}
Land, A.H., Doig, A.G.: An automatic method of solving discrete programming
  problems. Econometrica  \textbf{28}(3),  497--520 (1960)

\bibitem{Liberti2008}
Liberti, L.: Automatic generation of symmetry-breaking constraints. In:
  Combinatorial optimization and applications, Lecture Notes in Comput. Sci.,
  vol.~5165, pp. 328--338. Springer, Berlin (2008).
  \doi{10.1007/978-3-540-85097-7_31}

\bibitem{Liberti2012a}
Liberti, L.: Reformulations in mathematical programming: automatic symmetry
  detection and exploitation. Mathematical Programming  \textbf{131}(1-2),
  273--304 (2012). \doi{10.1007/s10107-010-0351-0},
  \url{http://dx.doi.org/10.1007/s10107-010-0351-0}

\bibitem{Liberti2012}
Liberti, L.: Symmetry in mathematical programming. In: Lee, J., Leyffer, S.
  (eds.) Mixed Integer Nonlinear Programming, IMA Series, vol.~154, pp.
  236--286. Springer, New York (2012)

\bibitem{LibertiOstrowski2014}
Liberti, L., Ostrowski, J.: Stabilizer-based symmetry breaking constraints for
  mathematical programs. Journal of Global Optimization  \textbf{60},  183--194
  (2014)

\bibitem{Margot2002}
Margot, F.: Pruning by isomorphism in branch-and-cut. Mathematical Programming
  \textbf{94}(1),  71--90 (2002). \doi{10.1007/s10107-002-0358-2},
  \url{http://dx.doi.org/10.1007/s10107-002-0358-2}

\bibitem{margot2003exploiting}
Margot, F.: Exploiting orbits in symmetric {ILP}. Mathematical Programming
  \textbf{98}(1),  3--21 (2003). \doi{10.1007/s10107-003-0394-6}

\bibitem{margot2010symmetry}
Margot, F.: Symmetry in Integer Linear Programming, chap.~17, pp. 647--686.
  Springer Berlin Heidelberg (2010). \doi{10.1007/978-3-540-68279-0_17}

\bibitem{MARTELLO1987213}
Martello, S., Toth, P.: Algorithms for knapsack problems. In: Martello, S.,
  Laporte, G., Minoux, M., Ribeiro, C. (eds.) Surveys in Combinatorial
  Optimization, North-Holland Mathematics Studies, vol.~132, pp. 213--257.
  North-Holland (1987). \doi{10.1016/S0304-0208(08)73237-7}

\bibitem{ostrowski2009symmetry}
Ostrowski, J.: Symmetry in integer programming. {PhD} dissertation, Lehigh
  University (2009)

\bibitem{OstrowskiEtAl2011}
Ostrowski, J., Linderoth, J., Rossi, F., Smriglio, S.: Orbital branching.
  Mathematical Programming  \textbf{126}(1),  147--178 (2011).
  \doi{10.1007/s10107-009-0273-x}

\bibitem{pfetsch2019computational}
Pfetsch, M.E., Rehn, T.: A computational comparison of symmetry handling
  methods for mixed integer programs. Mathematical Programming Computation
  \textbf{11}(1),  37--93 (2019). \doi{10.1007/s12532-018-0140-y}

\bibitem{PISINGER1999528}
Pisinger, D.: An exact algorithm for large multiple knapsack problems. European
  Journal of Operational Research  \textbf{114}(3),  528--541 (1999).
  \doi{10.1016/S0377-2217(98)00120-9}

\bibitem{salvagnin2005dominance}
Salvagnin, D.: A dominance procedure for integer programming. Master’s
  thesis, University of Padova, Padova, Italy  (2005)

\bibitem{Salvagnin2018}
Salvagnin, D.: Symmetry breaking inequalities from the {Schreier-Sims} table.
  In: van Hoeve, W.J. (ed.) Integration of Constraint Programming, Artificial
  Intelligence, and Operations Research. pp. 521--529. Springer International
  Publishing, Cham (2018)

\end{thebibliography}
\end{document}